\documentclass[preprint,12pt]{elsarticle_mod}

\usepackage{amsmath}
\usepackage{amssymb}
\usepackage{amsthm}
\usepackage{enumerate}

\usepackage{color}
\usepackage{graphicx}

\newtheorem{Theorem}[equation]{Theorem}

\newtheorem{Remark}[equation]{Remark}
\usepackage{setspace}
\numberwithin{equation}{section}
\numberwithin{table}{section}

\usepackage{graphicx}
 \usepackage{epsfig}

\usepackage{amssymb}
\usepackage{epsfig}
\usepackage{setspace}
\usepackage{amsthm}
\usepackage{amssymb}

\newcommand{\concunit}{g kg$^{-1}$ soil}

\begin{document}

\begin{frontmatter}



\title{Is adding charcoal to soil a good method for CO$_2$ sequestration? -- Modeling a spatially homogeneous soil}


\author[label1]{D.~Bourne}
\author[label2]{T.~Fatima}
\author[label2]{P.~van Meurs}
\author[label2,label3]{A.~Muntean}

\address[label1]{School of Mathematics and Statistics, University of Glasgow, UK}
\address[label2]{Department of Mathematics and Computer Science, Eindhoven University of Technology, The Netherlands}
\address[label3]{Institute for Complex Molecular Systems, Eindhoven University of Technology, The Netherlands}

\begin{abstract}
Carbon sequestration is the process of capture and long-term storage of atmospheric carbon dioxide ($CO_2$) with the aim to avoid dangerous climate change.  In this paper, we propose a simple mathematical model (a coupled system of nonlinear ODEs) to capture some of the dynamical effects produced by adding charcoal to fertile soils. The main goal is to understand to which extent charcoal is able to lock up carbon in soils. Our results are preliminary in the sense that we do not solve the $CO_2$ sequestration problem. Instead, we do set up a flexible modeling framework in which the interaction between charcoal and soil  can be tackled by means of mathematical tools.

We show that our model is well-posed and has interesting large-time behaviour. Depending on the reference parameter range (e.g. type of soil) and chosen time scale, numerical simulations suggest that adding charcoal typically postpones the release of $CO_2$.
\end{abstract}

\begin{keyword} Modeling  chemical kinetics in fertile soils\sep  Solvability of a nonlinear ODE system \sep Equilibria and steady states \sep Simulation \sep Biochar \sep  $CO_2$ sequestration


\end{keyword}

\end{frontmatter}


\section{Introduction}


In his Nature paper \cite{nature}, J. Lehmann argues that locking carbon up in soil makes more sense than storing it in plants and trees that eventually decompose, but does this idea work on a large timescale? A large community of soil scientists supports such ideas and attempts with experimental means to explore the sustainability  of adding charcoal (biochar) to soils; see for instance \cite{Glaser,Kolb, Steinbeiss,Celia} and
see also the review paper \cite{biota}. For more information on this research directions, often called the {\em Biochar project\footnote{{\em Biochar} := The idea of trapping carbon in soil for longer by storing it in the form of charcoal.}}, we refer the reader also to the sites \texttt{www.biochar-international.org} and \texttt{http://en.wikipedia.org/wiki/Biochar}. Briefly speaking, the Biochar project  aims at bringing clear advantages\footnote{Note also  the additional advantage of producing energy by burning organic matter to make charcoal.} (e.g. reduces soil greenhouse gas emissions, improves water and nutrient holding capacities, does not alter the carbon/nitrogen ratio, reduces soil acidity, removes pollutants), {\em but} is it a secure permanent solution?
 What about the possible negative effects like charcoal increases soil fertility and so increases the microbe population, which finally
 releases potentially more $CO_2$?
 It seems that there is no general agreement on whether putting charcoal in soil is generally a good idea or not. Therefore our interest.

In mathematical terms, our main question is:

{\em What is the large time behavior of the $CO_2$ dynamical system provoked by adding charcoal? }

In this context, the major issue is the complexity of the situation -- it is {\em a priori} not clear what would be a good "charcoal model" and to which extent capturing the effect of charcoal on $CO_2$ emissions is actually possible. This is the place where we wish to contribute.

\begin{figure}[h!]
\begin{centering}
                 \begin{tabular}{lcr}
                 \includegraphics[width=180pt]{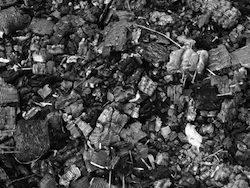}
                  \includegraphics[width=180pt]{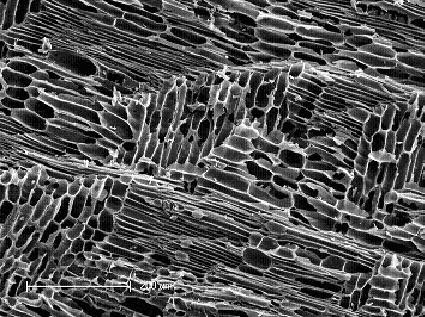}\\
                 \end{tabular}
                 \caption{Multiscale geometry of biochar (left: macro, right: micro). This is the place where nutrients, phenolics etc. undergo adsorption and desorption. }\label{fig:multiscale}
                 \end{centering}
\end{figure}

It is worth noting that charcoal is characterized by a very special porous  structure (see Figure \ref{fig:multiscale}), which is responsible for the high retention of water, dissolved organic nutrients, and even of pollutants  such as hydrocarbons and pesticides.
On top of this, the chemistry of soils is rather complex and precise (microscopic) characterizations of the microbial evolution are not available.  Furthermore, describing the transport of water together with nutrients, phenolics, pollutants (etc.) requires a good understanding of the heterogeneities of the soils.
Within this framework we treat a spatially homogeneous soil.  Herewith we avoid the aforementioned complications and propose a simple mathematical model, which is able to capture dynamical effects produced by adding charcoal to fertile soils. The model is  a nonlinearly coupled system of  deterministic ODEs which behaves well mathematically, that is  the system is solvable and its positive and bounded solution has a non-trivial large time behavior. Our main task is to explore the parameter space to investigate to which extent the presence of charcoal in soil affects $CO_2$ emissions on different time scales.


The paper is organized as follows: In Section \ref{CR} we describe mathematically chemical reactions in homogeneous media (here: soils) and propose a first model based on differential equations. We prove in Section \ref{analysis} that our model is well-posed in the sense of Hadamard and perform a stability analysis of the physically-relevant steady states. We illustrate the behavior of the profiles of the active concentrations and parameter effects in Section \ref{sim}. The effects observed regarding the addition of charcoal to soils are summarized in Section \ref{Discussion}. \ref{SAsub1} contains a discussion of the equilibria and stability of a reaction sub-block, while \ref{sens} reports on the sensitivity of $CO_2$ emission based on one of the reference parameter sets (Parameter Set 1).

We hope that our paper will bring the attention of the mathematical modeling community on the biochar issue. Note that, cf. Section \ref{Discussion}, there are many open modeling  aspects that would deserve a careful multi-disciplinary attention.

\section{Modeling chemical reactions in homogeneous fertile soils}
\label{CR}
\subsection{What happens if charcoal is added to soil?}

In this section we provide a simple model for the chemical reactions taking place in charcoal-enriched soil.
We model only those processes that are relevant to carbon dioxide emission: the break down of soil organic matter and charcoal by microbes and the subsequent release of carbon dioxide, the reproduction and death of the microbes, and the effect of charcoal on soil fertility.

We denote the species appearing in the chemical reactions by
\begin{equation}
\begin{array}{cl}
CO_2 & \textrm{carbon dioxide}, \\
Ch & \textrm{charcoal (actificially added to the soil)}, \\
Om & \textrm{soil organic matter (natural soil carbon)}, \\
M & \textrm{microbes}. \\
\end{array}
\end{equation}
Note that we do not distinguish between different types of soil organic matter (litter, recalcitrant organic matter, humus, etc.).
Also we only consider heterotrophic microbes, i.e., those that use organic carbon for growth.

Microbes in the soil break down the organic matter and charcoal (this is called
    mineralization), releasing the carbon, which then combines with oxygen to form carbon
    dioxide. Experimental evidence indicates that generally there is no shortage of oxygen in the soil. Having this mind we assume that oxygen is present everywhere in equal amounts and thus it enters our model as a parameter. We model the complex system of mineralization processes by means of the following chemical reactions mechanism:
\begin{align}
\label{CR1}
Om  &  \xrightarrow{k_1} n CO_2, \\
\label{CR2}
Ch & \xrightarrow{k_2} CO_2,
\end{align}
where $n>0$ is taken as a constant.
The reaction ``constants" $k_1$ and $k_2$ depend generally on the concentration of microbes, i.e,
$$k_i = k_i(M).$$
Here we assume that, as functions,  these reaction constants increase if the concentration of microbes increases.
Note that, in general, the reaction constants can also depend on other effects (like the concentration of phenolics in the soil), but for the sake of keeping things simple
we do not include these in our model.

The microbes need organic matter and oxygen to reproduce. Since we assumed that there is an abundance of oxygen, we can model the reproduction of microbes by means of
\begin{align}
\label{CR8}
M + \delta Om & \xrightarrow{k_3} (\mu + 1) M,
\end{align}
where $\delta, \mu >0$ are constants.
In general the reaction constant $k_3$ might depend on the fertility of the soil, which in turn
depends on the amount of charcoal in the soil.
For our theoretical investigations,  we neglect the intermediate step and assume directly that $k_3$ depends on the amount of charcoal,  $k_3 = k_3(Ch)$,
and that $k_3$ increases with charcoal concentration. However, note that the fertility of the soil contains so much {\em in situ } information that it  cannot be neglected in the practical design of a $CO_2$ sequestration scenario or if one wants to understand why {\em terra preta} (or `black earth' ) is so fertile.
Furthermore, in practice $k_3$ depends on many other factors, e.g., temperature, moisture, soil type, but we assume that these are all constant and
so they do not appear explicitly in our model.

We model the death of microbes by the chemical reaction
\begin{align}
M & \xrightarrow{k_4} \eta Om,
\label{CR10}
\end{align}
where $\eta > 0$ is a constant.

\subsection{Basics of chemical kinetics}

We denote the concentration of species $\mathcal{A}$ at time $t$ by $[\mathcal{A}](t)$, e.g., $[CO_2](t)$ is the
concentration of $CO_2$ in the soil at time $t$. In order to derive evolution equations for the
species concentrations we use the {\em simple reaction ansatz}, see, e.g., \cite{Atkins}. This assumption essentially states
that if our set of reactions is given by the mechanism
\begin{equation} \label{CK2}
    \sum_{i=1}^n \alpha_{ij} \mathcal{A}_i \xrightarrow{kj} \sum_{i=1}^n \beta_{ij} \mathcal{A}_i, \hspace{5mm} j = 1,
    \ldots, m,
\end{equation}

\noindent where $n \in \mathbb{N}$ denotes the number of species $A_i$, $m \in \mathbb{N}$ denotes
the number of chemical reactions, and $\alpha_{ij}, \beta_{ij} \in \mathbb{R}_+$ are stoichiometric coefficients, $k_j$ reaction constants,
then the {\em elementary reaction rates} are given by
\begin{equation} \label{CK4}
    r_j\left(\mathcal{A}_1,\mathcal{A}_2,\dots,\mathcal{A}_n\right) := k_j \prod_{i=1}^n [\mathcal{A}_i]^{ \alpha_{ij} }.
\end{equation}
Balancing the mass of the active species $\mathcal{A}_i$, we easily derive the evolution equations for the concentrations $[\mathcal{A}_i]$, viz.
\begin{equation} \label{CK3}
    \frac{d}{dt} [\mathcal{A}_i] = \sum_{j=1}^m (\beta_{i j} - \alpha_{i j}) r_j\left(\mathcal{A}_1,\mathcal{A}_2,\dots,\mathcal{A}_n\right), \hspace{5mm} i =
    1, \ldots, n.
\end{equation}

Before applying this methodology to \eqref{CR1}--\eqref{CR10}, we introduce a new notation, see Table \ref{table notation ui}, which is more convenient for the  analysis. For the sake of readability and clarity, we  use both notations throughout this paper.

\begin{table}[h!]
\centering
\begin{tabular}{|c|c|}
  \hline
  $u_1$ & $[Om]$ \\\hline
  $u_2$ & $[M]$ \\\hline
  $u_3$ & $[Ch]$ \\\hline
  $u_4$ & $[CO_2]$ \\
  \hline
\end{tabular}
\vspace{-2mm}
\caption{\footnotesize Alternative notation for the active concentrations.}
\label{table notation ui}
\end{table}
\begin{Remark}(Restriction to spatially homogeneous soils)
Within the framework of this paper, we consider a ``continuously stirred tank reactor" case, a scenario intensively used in chemical engineering; see, e.g., \cite{Atkins}.  In terms of soils, this means that we focus our modeling on a single space location, where the measurements are  made, and we follow how the information ``flows" over  physically-important timescales.  To this end, we assume the soil to be homogeneous in the sense that no spatial substructures (typically called microstructures) appear, i.e., all soil components (gravel, sand, solid nutrients, water, etc) are well-mixed.   We postpone for later the study of the more realistic case when the soil heterogeneities will be explicitly taken into account in terms of porosities, tortuosities, permeabilities very much in the spirit of  \cite{Bear} (general theory of flows in porous media), \cite{Tasnim, Ekeoma, Mariya} (multiscale approaches to the chemical corrosion of concrete, smoldering combustion and plant growth, respectively), \cite{Verma} (accumulation of cadmium in plants). Also, at a later stage it would be interesting to study the effect of the charcoal's platelet-like microstructure (see Figure \ref{fig:multiscale}) on the efficiency of adsorption and desorption of the nutrients. Most likely this would lead to a two-scale ODE system intimately coupled with evolution equations for the transport and storage of nutrients.
\end{Remark}

Applying the {\em simple reaction ansatz} to \eqref{CR1}--\eqref{CR10}, and assuming additionally that the system has a constant source $s \ge0$ of organic matter, yields the nonlinear coupled system of ODEs
\begin{align}
\label{DE8}
\frac{d}{dt} u_1 & = -k_1 (u_2) u_1 - \delta k_3 (u_3) u_2 u_1^\delta + \eta k_4 u_2 + s, \\
\label{DE9}
\frac{d}{dt} u_2 & =  \mu k_3 (u_3) u_2 u_1^\delta - k_4 u_2,\\
\label{DE3}
\frac{d}{dt} u_3 & = -k_2 (u_2) u_3, \\
\label{DE6}
\frac{d}{dt} u_4 & = n k_1 (u_2) u_1 + k_2 (u_2) u_3.
\end{align}
The source $s$ can be thought of as organic matter entering the soil from the surface in the form of dead leaves, plants, etc. This system also requires initial conditions. Their role is to incorporate the type of soil.
Throughout the rest of this paper we study the system \eqref{DE8}--\eqref{DE6}.


\section{Mathematical analysis of the system \eqref{DE8}--\eqref{DE6}}\label{analysis}

We start by introducing a set of assumptions on the model parameters entering \eqref{DE8}--\eqref{DE6}.
These assumptions will be used to prove global existence of positive and bounded concentrations $u_i$ and to study the steady states of this nonlinear ODE system.

\subsection{Restrictions on the model parameters}\label{restrictions}

We assume that
\begin{equation}\label{PR05}
\delta \geq 1.
\end{equation}
Assumption \eqref{PR05}, together with the assumptions given below on the constitutive functions $k_i$, ensure
that the right-hand side of the system \eqref{DE8}--\eqref{DE6} is Lipschitz continuous, which guarantees that our ODE system admits a unique local classical solution.

In addition to choosing that $\delta, \eta, \mu, n > 0$, we also assume that
\begin{equation} \label{dgem}
    \delta \geq \eta \mu.
\end{equation}
The physical meaning of (\ref{dgem}) is explained in \ref{SAsub1}.
The condition (\ref{dgem}) is used in Section \ref{LB} to ensure that the solution to \eqref{DE8}--\eqref{DE6} does not blow-up in finite time.

Since the $k_i$ are reaction constants, we assume that they satisfy $k_i > 0$ for $i\in\{1,2,3,4\}$. Note however that $k_i$ are nearly never true constants; they often incorporate a certain dependence on important physical/environmental quantities (here:  spatial location, temperature, soil fertility, oxygen content, water content, etc). Here we take $k_4$ to be  constant and assume that the functions $k_i : \mathbb{R} \rightarrow (0,\infty)$, $i \in \{1,2,3\}$, are Lipschitz continuous and strictly increasing. For example, $k_1$ being strictly increasing means that an increase of microbes in the soil leads to an increase in the rate of break down of organic matter.

Finally, we assume that the initial concentrations are positive and bounded, i.e. $u_i(0) = u_i^0 \in [0,\infty)$, $i \in \{1,2,3,4\}$.

\subsection{Positivity of concentrations}
\label{P}

In this section we show that the concentrations $u_1,u_2,u_3,u_4$ are nonnegative for all times if their initial values are nonnegative. It suffices to show for each $i\in\{1,2,3,4\}$ that if $u_i=0$ and $u_j \ge 0$ for all $j \ne i$, then $\dot{u}_i\ge 0$. This turns to be a trivial exercise:
\begin{align*}
\dot{u}_1(0,u_2,u_3,u_4) & = \eta k_4 u_2 + s \geq 0, \\
\dot{u}_2(u_1,0,u_3,u_4) & = 0,\\
\dot{u}_3(u_1,u_2,0,u_4) & = 0, \\
\dot{u}_4(u_1,u_2,u_3,0) & = n k_1 (u_2) u_1 + k_2 (u_2) u_3 \ge 0.
\end{align*}

\subsection{$L^\infty$ bounds on concentrations}
\label{LB}
We prove that the concentrations $u_i$ do not blow-up in finite time.
Fix arbitrary initial conditions $u_i^0$. Then, based on the result of Section \ref{P}, we can assume that $u_i \ge 0$ for all $i=1,2,3,4$.

From the positivity of $u_i$ and $k_i$, it follows immediately from \eqref{DE3} that
\begin{equation}\label{LB1}
    \| u_3 \|_\infty \leq u_3^0.
\end{equation}
Adding equation \eqref{DE8} to $\eta$ times equation \eqref{DE9} gives
\begin{equation}\label{LB2}
\begin{aligned}
    \frac{d}{dt} (u_1 + \eta u_2)
    = -k_1 (u_2) u_1 - (\delta - \eta \mu) k_3 (u_3) u_2 u_1^\delta +s \leq s.
\end{aligned}
\end{equation}
The inequality \eqref{LB2} follows from \eqref{dgem} and the positivity of the $k_i$ and $u_i$. From \eqref{LB2} we conclude that $u_1$ and $u_2$ satisfy $L^\infty$ bounds on any finite time interval. The numerics suggest that this bound is independent of the length of this time interval, but we do not need this here; see section \ref{sim}.

Relying on the $L^\infty$ bounds on $u_i$ for $i \in \{ 1,2,3 \}$ on any finite time interval $[0,\tau]$, we can bound the right-hand side of \eqref{DE6} by a constant $C(\tau)$. 
Integration yields the bound
\begin{equation}\label{LB5}
    u_4(t) \leq C(\tau) t + u_4^0
\end{equation}
for all $t \in [0,\tau]$,
which immediately gives a bound on $u_4$ on any time interval $[0,\tau]$.

\subsection{Well-posedness}
Based on the positivity and the $L^\infty$ bounds on concentrations, together with the Lipschitz continuity of the right-hand side of \eqref{DE8}--\eqref{DE6}, we recall classical ODE theory  (see \cite{Coddington,Amann}, e.g.) to prove the following result:

\begin{Theorem} \label{thm WP}
(Global solvability). Assume that the assumptions stated in section \ref{restrictions} hold. Then
for any set of initial conditions $u_i(0)=u_i^0 \ge 0$, 
the system \eqref{DE8}--\eqref{DE6} has a unique classical solution $u_i : [0,\infty) \rightarrow \mathbb{R}$, $i \in \{ 1,2,3,4 \}$.
\end{Theorem}

Furthermore, a Gronwall-like argument can be employed to show that this classical solution depends continuously on the initial data and {\em all} model parameters. Since this argument is rather standard, we omit to show it here.

\subsection{Equilibria and stability of the  system \eqref{DE8}--\eqref{DE3}}\label{NSSA}

First note that $u_4$ does not appear in the right hand side of \eqref{DE8}--\eqref{DE6}. Hence equation \eqref{DE6} decouples from the system, in the sense that we do not need \eqref{DE6} to solve the subsystem \eqref{DE8}--\eqref{DE3}. Having this in mind, it is sufficient to  study the equilibria of the reduced system \eqref{DE8}--\eqref{DE3}. The reader is referred to \ref{SAsub1}
for a discussion of the equilibria and stability of the reaction block given by \eqref{CR8} and \eqref{CR10}. For basic notions of dynamical systems\footnote{Dynamical systems theory proved to be very successful in a series of cases arising in biology and ecology; compare  for instance \cite{Feng,Rui} and references cited therein. We expect therefore that  dynamical systems delivers results in the case of  biochar research as well.}, see \cite{Guck}, e.g.


We first search for the equilibria of the decoupled system given by \eqref{DE8}, \eqref{DE9} and \eqref{DE3}. By equating the right-hand side of \eqref{DE3} to zero, it follows that $u_3 = 0$. By substituting this into equations \eqref{DE8} and \eqref{DE9}, we obtain
\begin{align}\label{NS4}
    0 &= -k_1 (u_2) u_1 - \delta k_3 (0) u_2 u_1^\delta + \eta k_4 u_2 + s, \\\label{NS5}
    0 &= \big( \mu k_3 (0) u_1^\delta - k_4 \big) u_2.
\end{align}
For convenience we write $k_3$ instead of $k_3(0)$ in the remainder of this section. Equation \eqref{NS5} is satisfied if and only if either
\begin{gather}\label{NS6}
    u_2 = 0, \textrm{ or} \\\label{NS7}
    u_1 = \bigg( \frac{ k_4 }{ \mu k_3 } \bigg)^\frac1\delta =: U_1.
\end{gather}
Let us treat the two cases separately:

\begin{description}
  \item[Case \eqref{NS6}:] It immediately follows from \eqref{NS4} that $u_1 = s/k_1(0)$.
  \item[Case \eqref{NS7}:] By inserting \eqref{NS7} in \eqref{NS4} we get
  \begin{equation} \label{NS8}
    0 = -k_1 (u_2) U_1 - \frac{k_4}\mu (\delta - \eta \mu) u_2 + s.
  \end{equation}

  \noindent The right-hand side of \eqref{NS8} is strictly decreasing as a function of $u_2$. Hence it has at most one solution $u_2$. A necessary condition for the existence of such a solution is that the right-hand side is nonnegative for $u_2 = 0$. This is the case when
  \begin{equation}\label{NS9}
    s \geq k_1(0) U_1.
  \end{equation}
  From now on we assume that the $k_i$ and the parameters $\delta, \eta, \mu$ are chosen such that \eqref{NS8} has a solution whenever \eqref{NS9} holds. We will call this solution $u_2^\ast$. For example, a solution exists if \eqref{NS9} holds and $\delta - \eta \mu > 0$.
\end{description}

Therefore, depending on the parameter $s$, we have one or two equilibrium points: If $s \leq k_1(0) U_1$, then we have only one equilibrium point $u_e^1$ given by
\begin{equation}\label{NS10}
    u_e^1 := (u_1, u_2, u_3) = \bigg( \frac{s}{k_1(0)}, 0, 0 \bigg).
\end{equation}

If $s > k_1(0) U_1$, we have the additional equilibrium point $u_e^2$ given by
\begin{equation}\label{NS11}
    u_e^2 := (u_1, u_2, u_3) = \bigg( U_1, u_2^\ast, 0 \bigg),
\end{equation}
where $u_2^\ast$ satisfies \eqref{NS8}. Therefore $s = k_1(0) U_1$ is a bifurcation point.

To test the stability of the equilibrium points $u_e^1$ and $u_e^2$, we linearize the system \eqref{DE8}, \eqref{DE9}, \eqref{DE3}. Let $J$ denote the Jacobian matrix of this system. A brief calculation shows that
\begin{equation}\label{NS12}
    J(u_e^1) = \begin{bmatrix}
        -k_1(0) & -\frac{k_1'(0)}{k_1(0)} s
                  - \delta k_3 \big( \frac{s}{k_1(0)} \big)^\delta
                  + \eta k_4
                    & 0       \\
        0       & \mu k_3 \big( \frac{s}{k_1(0)} \big)^\delta - k_4
                    & 0       \\
        0       & 0 & -k_2(0) \\
      \end{bmatrix}.
\end{equation}

\noindent The eigenvalues of $J(u_e^1)$ are given by the entries on the diagonal. The eigenvalues $-k_1(0)$ and $-k_2(0)$ are negative, whereas the sign of the third eigenvalue changes from negative to positive as $s$ passes the bifurcation point. So $u_e^1$ is asymptotically stable if $s < -k_1(0) U_1$ and is unstable if $s > -k_1(0) U_1$.

We follow the same procedure for $u_e^2$. First we obtain
\begin{equation}\label{NS13}
    J(u_e^2) = \begin{bmatrix}
        -k_1(u_2^\ast) - \delta^2 k_3 u_2^\ast U_1^{\delta-1}
        & -k_1'(u_2^\ast) U_1 - \frac{k_4}{\mu} (\delta - \eta \mu)
          & - \frac\delta\mu \frac{k_3'(0)}{k_3(0)} k_4 u_2^\ast \\
        \delta \mu k_3 u_2^\ast U_1^{\delta-1}
        & 0
          & \frac{k_3'(0)}{k_3(0)} k_4 u_2^\ast \\
        0
        & 0
          & -k_2(u_2^\ast) \\
      \end{bmatrix}.
\end{equation}
Let us denote the $2 \times 2$ upper-left block of $J(u_e^2)$ by
\begin{equation*}
    \begin{bmatrix}
        A_1 & A_2 \\
        A_3 & 0   \\
      \end{bmatrix}.
\end{equation*}
Note that $A_1, A_2 < 0$ and $A_3 > 0$. Therefore the eigenvalues of $J(u_e^2)$ are
%
\begin{equation*}
    -k_2(u_2^\ast), \hspace{3mm} \frac{A_1}2 + \sqrt{ \frac{A_1^2}4 + A_2 A_3 }, \hspace{3mm} \textrm{and } \frac{A_1}2 - \sqrt{ \frac{A_1^2}4 + A_2 A_3 }.
\end{equation*}
Since $A_1 < 0$ and $A_2 A_3 < 0$, the real parts of all the three eigenvalues are negative, which proves that $u_e^2$ is asymptotically stable.

In summary, for each $s > 0$
there is one stable equilibrium of the decoupled system \eqref{DE8}, \eqref{DE9}, \eqref{DE3}. Depending on the size of the source $s$, this equilibrium is either given by \eqref{NS10} or by \eqref{NS11}.
Note that the full system \eqref{DE8}--\eqref{DE6} does not have any equilibrium points since $\dot{u_4}>0$ (unless $s=0$, in which case $u_i = 0$ for all $i=1,2,3,4$ is an equilibrium).

\section{Nondimensionalisation}
\label{N}
Before solving the system numerically, we rescale it (very much in the spirit of \cite{Segel}).
We consider the following scalings for the time, concentrations, and reaction rates:
$t=\tau \tilde{t}$, where $\tau$ is the reference time,
$u_i=U_i \tilde{u}_i$, where $U_i$ is the reference concentration of species $i$, and $k_i = K_i \tilde{k}_i$, where $K_i$ is the reference reaction constant.
Substituting these into equations \eqref{DE8}--\eqref{DE6} gives
\begin{equation}
\label{N1}
\begin{aligned}
\frac{d}{d \tilde{t}}  \tilde{u}_1 & = -\tau_1\tilde{k_1}\tilde{u}_1- \tau_2
\tilde{k_3}\tilde{u}^\delta_1\tilde{u}_2+ \tau_3\tilde{k_4}\tilde{u}_2 + \tau_4 s, \\
\frac{d}{d \tilde{t}}  \tilde{u}_2 & =  \tau_5
\tilde{k_3}\tilde{u}^\delta_1\tilde{u}_2-\tau_6\tilde{k_4}\tilde{u}_2, \\
\frac{d}{d \tilde{t}}  \tilde{u}_3 & =  -\tau_7\tilde{k_2}\tilde{u}_3, \\
\frac{d}{d \tilde{t}}  \tilde{u}_4 & =  \tau_8 \tilde{k_1}\tilde{u}_1+\tau_9
\tilde{k_2}\tilde{u}_3,
\end{aligned}
\end{equation}
where $\tau_\alpha, \alpha\in\{1,2,\dots ,9\}$, denote the characteristic time scales. Table \ref{tau} lists their dependence on the reference constants.

\begin{table*}[h!]
  \centering
  \begin{tabular}{|c|c|} \hline
  Characteristic  & Definition\\
  time scale& \\
   \hline
    $\tau_1$ & $\tau {K_1}$  \\
    \hline
      $\tau_2$ & $\tau \delta{K_3}U_1^{\delta-1} U_2$  \\
    \hline
    $\tau_3$ & $\tau \eta{K_4}U_2 U_1^{-1}$  \\
          \hline
    $\tau_4$ & $\tau U_1^{-1}$\\
        \hline
    $\tau_5$ & $\tau \mu {K_3}U_1^{\delta}$  \\
        \hline
    $\tau_6$ & $\tau {K_4}$  \\
        \hline
    $\tau_7$ & $\tau {K_2}$  \\
        \hline
    $\tau_8$ & $\tau n K_1 U_1 U_4^{-1}$  \\
        \hline
    $\tau_9$ & $\tau {K_2}U_3 U_4^{-1}$\\
    \hline
  \end{tabular}
  \caption{List of the involved characteristic time scales.}
  \label{tau}
\end{table*}

\section{Numerical simulation of the system \eqref{DE8}--\eqref{DE6}}\label{sim}
\subsection{The philosophy of our approach}
Here we illustrate numerically the behaviour of the solution to our ODE system. The main interest lies in predicting how the emission of $CO_2$ into the atmosphere changes if we put charcoal in the soil \cite{nature}. As the parameter space is large, we start with a typical parameter set for our model (i.e. as many parameters are $\mathcal{O}(1)$). After discussing these results, we test our model with physical parameters for $U_i$ and $\tau$, respectively the reference values for the concentrations and time.

We start by choosing the following linear constitutive functions for the reaction rates:
\begin{equation*}
\begin{aligned}
\label{ks}
\tilde{k}_1 (\tilde{u}_2) & := a_1 \tilde{u}_2 + b_1, \quad \tilde{k}_2 (\tilde{u}_2) : = a_2 \tilde{u}_2 + b_2, \\
\tilde{k}_3 (\tilde{u}_3) &:= a_3 \tilde{u}_3 + b_3, \quad \tilde{k}_4  := b_4.
\end{aligned}
\end{equation*}
Next, we wish that our model inherits the following features (ordering of characteristic time-scales):
\begin{itemize}
  \item Equality in \eqref{dgem} (i.e. the reproduction and death of microbes should be balanced).
  \item The reproduction and death of microbes is considerably faster than the break down of organic matter. The latter is still faster than the break down of charcoal.
  \item There should be a stable equilibrium with $\tilde{u}_1, \tilde{u}_2 > 0$ (see \eqref{NS11}).
\end{itemize} The three features introduce natural constraints in the parameter space.
Note that the first feature is satisfied if
\begin{equation} \label{NSN1}
    \delta = \eta \mu,
\end{equation}
while to fulfill the second one  we need at least
\begin{equation}\label{NSN2}
    K_2 < \frac1{2} K_1 < \frac1{2} \min \{K_3, K_4\}.
\end{equation}
Finally, to ensure that the equilibrium as given in \eqref{NS11} is stable, we take
\begin{equation}\label{NSN3}
    s = \alpha K_1 b_1 U_1,
\end{equation}
where $\alpha > 1$ is still to be determined. Further, we normalize the system such that $1 = \tilde{u}_1^e = \tilde{u}_2^e$. By definition of $U_1$ (see~\eqref{NS7}), $\tilde{u}_1^e = 1$ is automatically satisfied. As we like to have $U_1$ as a reference value that we can choose, we change~\eqref{NS7} into the following condition on $K_4$:
\begin{equation}\label{NSN4}
    K_4 := U_1^\delta K_3 \mu \frac{b3}{b4}.
\end{equation}
Last, we need
\begin{equation}\label{NSN5}
    1 = \tilde{u}_2^e = \frac{ s - K_1 U_1 b_1 }{ K_1 a_1 U_1 + K_4 b_4 U_2 (\delta - \eta \mu)/\mu } = (\alpha - 1)\frac{b_1}{a_1}.
\end{equation}

\begin{table*}
  \centering
  \begin{tabular}{|c|c||c|c|c|}\hline
     parameter & value & ref. constant & value & unit \\ \hline
     $a_1$ &  1 &    $K_1$ & 0.01 & s$^{-1}$ \\
     $a_2$ &  1 &    $K_2$ & $10^{-3}$ & s$^{-1}$ \\
     $a_3$ &  1.9 &  $K_3$ & 1 & m$^{3 \delta}$ mol$^{-\delta}$ s$^{-1}$\\
     $b_1$ & 1 &     $U_1$ & 1 & mol m$^{-3}$ \\
     $b_2$ & 1 &     $U_2$ & 1 & mol m$^{-3}$ \\
     $b_3$ & 0.1 &   $U_3$ & 1 & mol m$^{-3}$ \\
     $b_4$ & 1 &     $U_4$ & $10^3$ & mol m$^{-3}$ \\
     $\mu$ & 1 &     $\tau$ & 1 & s  \\
     $\delta$ & 10 &  & & \\
     $n$ & 10 &       & & \\\hline
  \end{tabular}
  \caption{Parameter Set 1.}
  \label{tab:params}
\end{table*}

The parameters and reference values that are still free for us to choose, are listed in Table \ref{tab:params}, together with the values we chose for them. With this set of values and \eqref{NSN1} -- \eqref{NSN5} we obtain $\eta = 10$, $s = 0.02$ mol m$^{-3}$ s$^{-1}$ and~$K_4 = 0.1$ s$^{-1}$.


The values in Table \ref{tab:params} were initially set to be $1$, except for $K_1$, $K_2$, $\delta$ and $n$. The reference constants $K_1$ and $K_2$ are chosen to satisfy \eqref{NSN2}. By taking $\delta = 10$, we model that microbes need to break down, on average, ten organic matter particles before they reproduce. We put $n = 10$ to model that ten $CO_2$ molecules are produced during the mineralization process of a single organic matter particle.

The resulting simulation showed minor response in the values for $\tilde{u}_i$ when charcoal was added, so we altered the value for the parameters $a_3$ and $b_3$ to make the reproduction of microbes more dependent on $[Ch]$, while conserving $\tilde{k}_2(1) = \tilde{k}_i(1) = 2$. The value for $U_4$ only effects the scaling for $\tilde{u}_4$; it is chosen such that it is $\mathcal{O}(1)$ on the long time scale (i.e. at which the system converges back to equilibrium).

\begin{figure}[h!]
\begin{centering}
                \includegraphics[width=350pt]{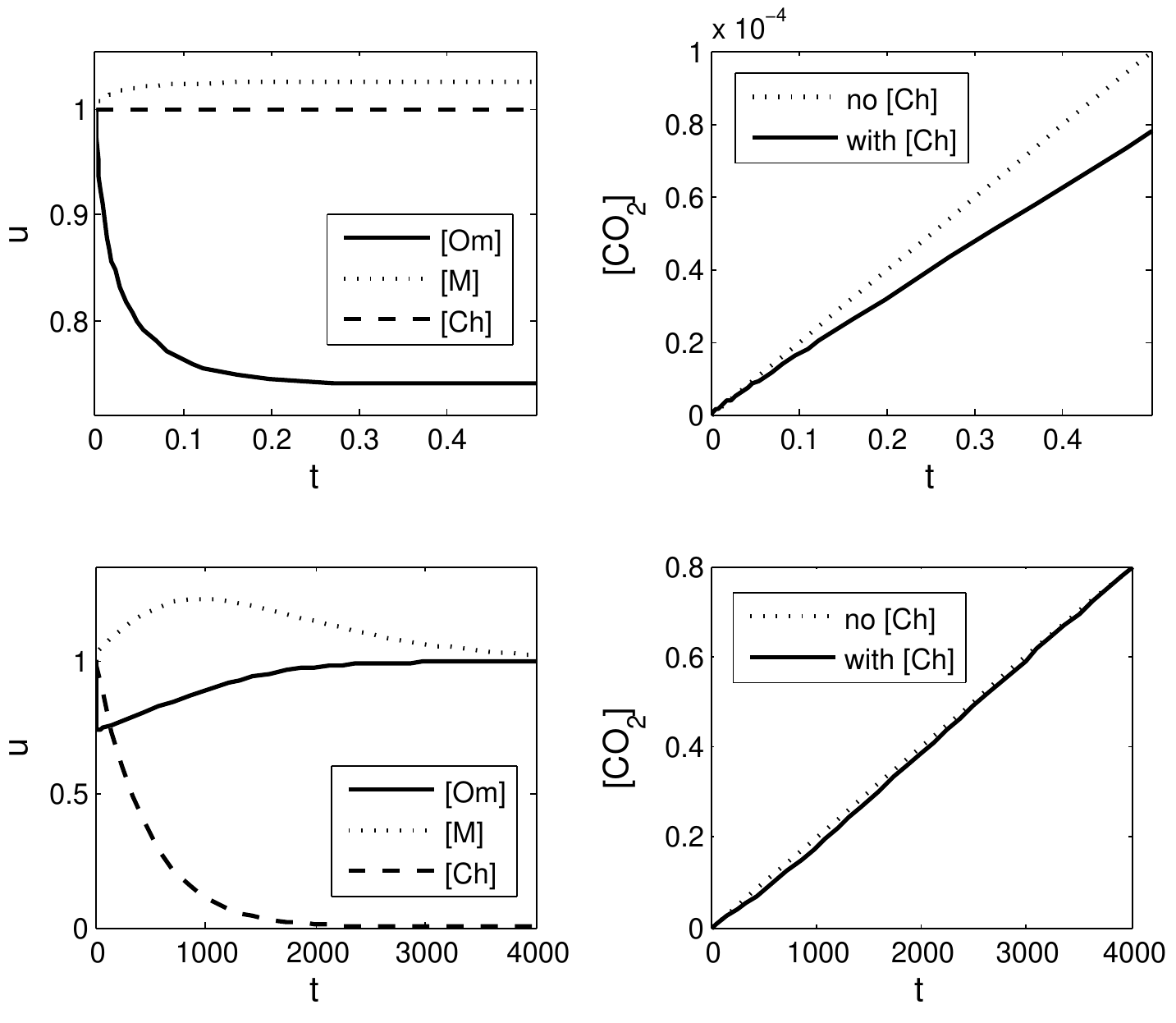}
                 \caption{These figures show the short-time (top) and long-time (bottom) behaviors of the system initially at equilibrium. Charcoal is added at time $\tilde{t} = 0$. The graph of $CO_2$ is put into another plot, together with the $CO_2$ emission that would occur if no charcoal was added to the soil.}\label{pmsetone}
                 \end{centering}
\end{figure}

Figure \ref{pmsetone} shows the result of the simulation with the parameters and reference values as in Table \ref{tab:params} (where the remaining parameters and reference values are computed via \eqref{NSN1} -- \eqref{NSN5}). The initial value is $\tilde{u}_1(0) = \tilde{u}_2(0) = \tilde{u}_3(0) = 1$ and $\tilde{u}_4(0) = 0$ (i.e. at time $\tilde{t} = 0$ charcoal is added to the soil in the otherwise stable state.) Although the simulation is carried out for the dimensionless $\tilde{u}_i$, we will refer to them by $[Om]$, $[M]$, $[Ch]$ and $[CO_2]$ for clarity. Figure \ref{pmsetone} shows various interesting phenomena:

\begin{itemize}
  \item The concentrations $[Om]$ and $[M]$ change on a short time scale ($\mathcal{O}(\tilde{t}) = 0.1$). Essentially this is because their time derivatives depend on $[Ch]$ through $k_3$.
  \item On an intermediate time scale $\mathcal{O}(\tilde{t}) = 100$, the $CO_2$ emission decreases when charcoal in put in the soil.
  \item On a long time scale ($\mathcal{O}(\tilde{t}) = 2000$), $[Ch]$ decreases exponentially fast to 0. Therefore $[Om]$ and $[M]$ converge back to their initial, equilibrium values.
  \item On the same long time scale, the $[CO_2]$ emission is almost the same as in the case in which no charcoal is added.
\end{itemize}

The graphs of the $[CO_2]$ emission can be explained by two effects (see fourth equation of \eqref{N1}). One term comes from the mineralization of charcoal with $CO_2$ as by-product; it increases the $[CO_2]$ emission. The other term comes from the mineralization of organic matter. So because $[Om]$ decreases if charcoal is added, this has a decreasing effect on the $[CO_2]$ emission. The ratio of $\tau_8$ and $\tau_9$ determines how much these two effects matter which respect to one other.

Still we like to understand the sensitivity of the $CO_2$ emission on the parameter space on a more detailed level. We elaborate on this further in \ref{sens}.

\subsection{Realistic parameters}

\begin{table*}
  \centering
  \begin{tabular}{|c|c||c|c|c|}\hline
     parameter & value & ref. constant & value & unit \\ \hline
     $a_1$ &  1 &    $K_1$ & $5 \cdot 10^{-8}$ & s$^{-1}$ \\
     $a_2$ &  1 &    $K_2$ & $2 \cdot 10^{-8}$ & s$^{-1}$ \\
     $a_3$ &  1.9 &  $K_3$ & $3 \cdot 10^{-10}$ & (\concunit )$^{-\delta}$ s$^{-1}$\\
     $b_1$ & 1 &     $U_1$  & $18$\cite{Steiner}  & \concunit \\
     $b_2$ & 1 &     $U_2$  & $0.2$  \cite{Steiner} & \concunit \\
     $b_3$ & 0.1 &   $U_3$ & $100$ \cite{Steiner}  & \concunit \\
     $b_4$ & 1 &     $U_4$ & $180$ \cite{Smith} & \concunit \\
     $\mu$ & 1 &     $\tau$ & $1$ & year  \\
     $\delta$ & 2 &  & & \\
     $n$ & 10 &       & & \\\hline
  \end{tabular}
  \caption{Parameter Set 2.}
  \label{tab:params:2}
\end{table*}

In contrast to putting as much parameters and reference values equal to $1$, we now take characteristic values for $U_i$ from literature (cf. e.g. \cite{Steiner,Smith}). Furthermore, we tune the time scale such that charcoal is broken down in the order of $1$ year. By exploring the parameter space in this way, we noticed that $[M]$ grows unnaturally fast (by a factor $10$ on the short time scale). By taking $\delta = 2$ and playing with the values for $K_i$, we could reduce it to a more physical growing factor. Furthermore, we increased the value for $U_4$ considerably (rather than using $U_4 = 1.26$ \concunit{} \cite{Smith}, as measured for a one-year period); this only changes the value of $\tilde{u}_4$ by a constant. The resulting set of parameters is shown in Table \ref{tab:params:2}. With this set of values and \eqref{NSN1} -- \eqref{NSN5} we obtain $\eta = 2$, $s = 1.8 \cdot 10^{-6}$ \concunit{} s$^{-1}$ and $K_4 = 9.72 \cdot 10^{-9}$ s$^{-1}$. Note that we also need to satisfy \eqref{NSN2}, where we now change units from molecular to mass concentrations.

\begin{figure}[h!]
\begin{centering}
                  \includegraphics[width=350pt]{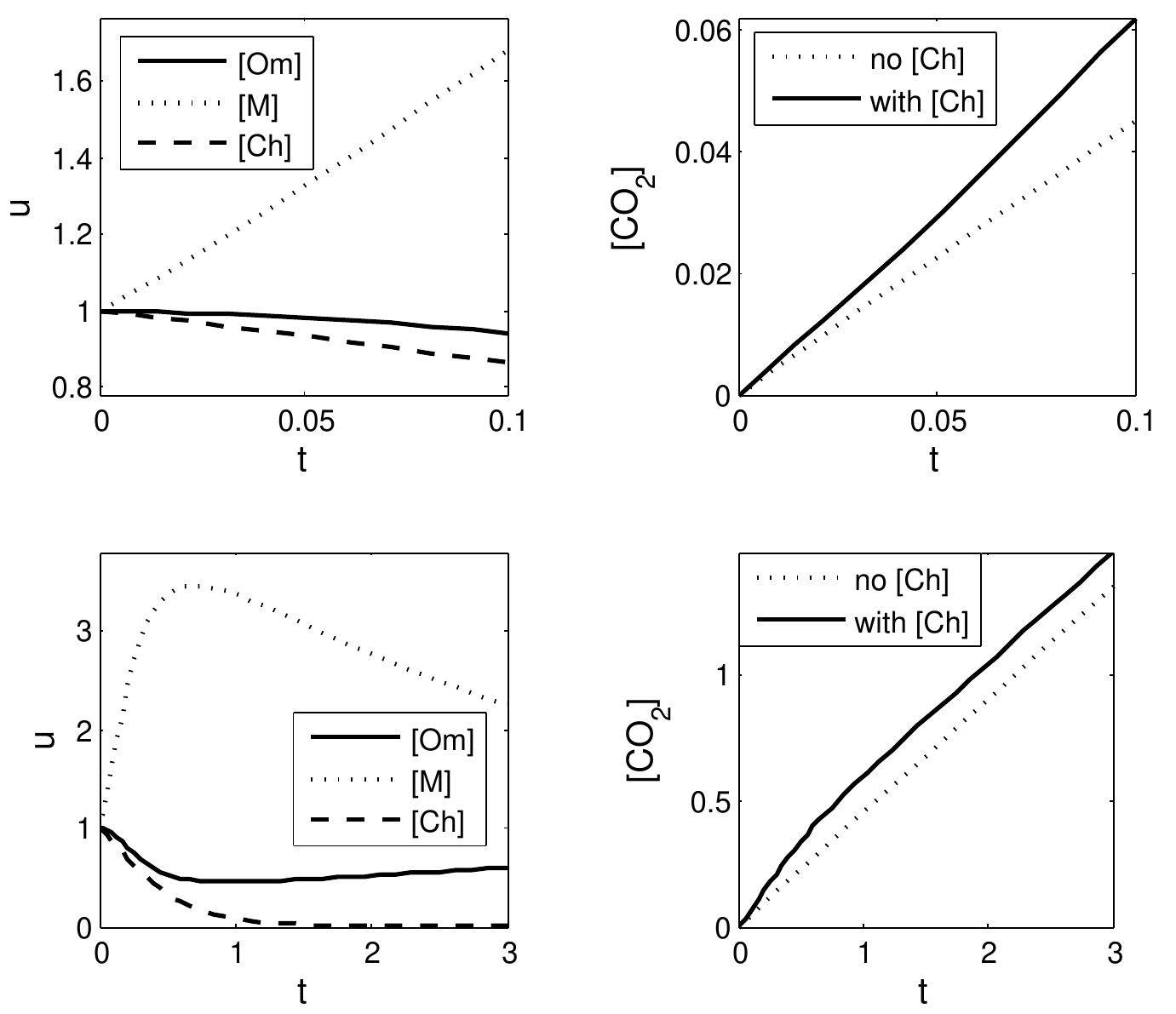} 
                 \caption{\label{pmsettwo}
                 Result of the simulation by using Parameter Set 2 (see Table \ref{tab:params:2}).}
                 \end{centering}
\end{figure}

Figure \ref{pmsettwo} shows the results from the simulation. In comparison to Figure \ref{pmsetone}, the most remarkable difference is that there is no response of $\tilde{u}_i$ on a short time scale, even though one would expect this. Furthermore, the $[CO_2]$ emission increases when charcoal is added. We kept on seeing this while exploring the parameter space. In the next section we try to connect these results to those of the previous parameter set.

\subsection{Further insight in the parameter space}

We aim to get similar results as depicted in Figure \ref{pmsetone} by deviating only a little from the parameter and reference values as in Table \ref{tab:params:2}.

We start by taking $\delta = 5$, which is more reasonable than $\delta = 2$ (see~\eqref{CR8}). From Table \ref{tau}, we see that $\tau_5$ scales with $U_1^\delta$, so it becomes large. We believe this is the reason for the non-physical increase in $[M]$. To prevent $\tau_5$ from being too large, we divide the previous value of $U_1$ by 5. This gives $U_1 = 3.6$ \concunit{}, which means that the soil contains less organic matter. We can further decrease $\tilde{u}_2$ by enlarging $U_2$. As a result, we take $U_2 = 2$ \concunit, i.e. we take a soil with 10 times as much microbes.

\begin{table*}[h!]
  \centering
  \begin{tabular}{|c|c||c|c|c|}\hline
     parameter & value & ref. constant & value & unit \\ \hline
     $a_1$ &  1 &    $K_1$ & $5 \cdot 10^{-8}$ & s$^{-1}$ \\
     $a_2$ &  1 &    $K_2$ & $5 \cdot 10^{-8}$ & s$^{-1}$ \\
     $a_3$ &  1.9 &  $K_3$ & $5 \cdot 10^{-8}$ & (\concunit )$^{-\delta}$ s$^{-1}$\\
     $b_1$ & 1 &     $U_1$ & 3.6 & \concunit \\
     $b_2$ & 1 &     $U_2$ & 2 & \concunit \\
     $b_3$ & 0.1 &   $U_3$ & 100 & \concunit \\
     $b_4$ & 1 &     $U_4$ & $180$ & \concunit \\
     $\mu$ & 1 &     $\tau$ & 1 & year  \\
     $\delta$ & 5 &  & & \\
     $n$ & 100 &       & & \\\hline
  \end{tabular}
  \caption{Parameter Set 3.} 
  \label{tab:params:3}
\end{table*}

By taking $n = 100$ (see \eqref{CR1} for the interpretation), we make the $CO_2$ emission more dependent on the mineralization of organic matter rather than the mineralization of charcoal. Finally, we tune the values for $K_i$ a bit to resemble the results as shown in Figure \ref{pmsetone}. Table \ref{tab:params:3} shows the list of parameters and reference values. Together with \eqref{NSN1} -- \eqref{NSN5} we obtain $\eta = 5$, $s = 3.6 \cdot 10^{-7}$ \concunit{} s$^{-1}$ and $K_4 = 3.0 \cdot 10^{-6}$ s$^{-1}$.

\begin{figure}[h!]
\begin{centering}
                  \includegraphics[width=350pt]{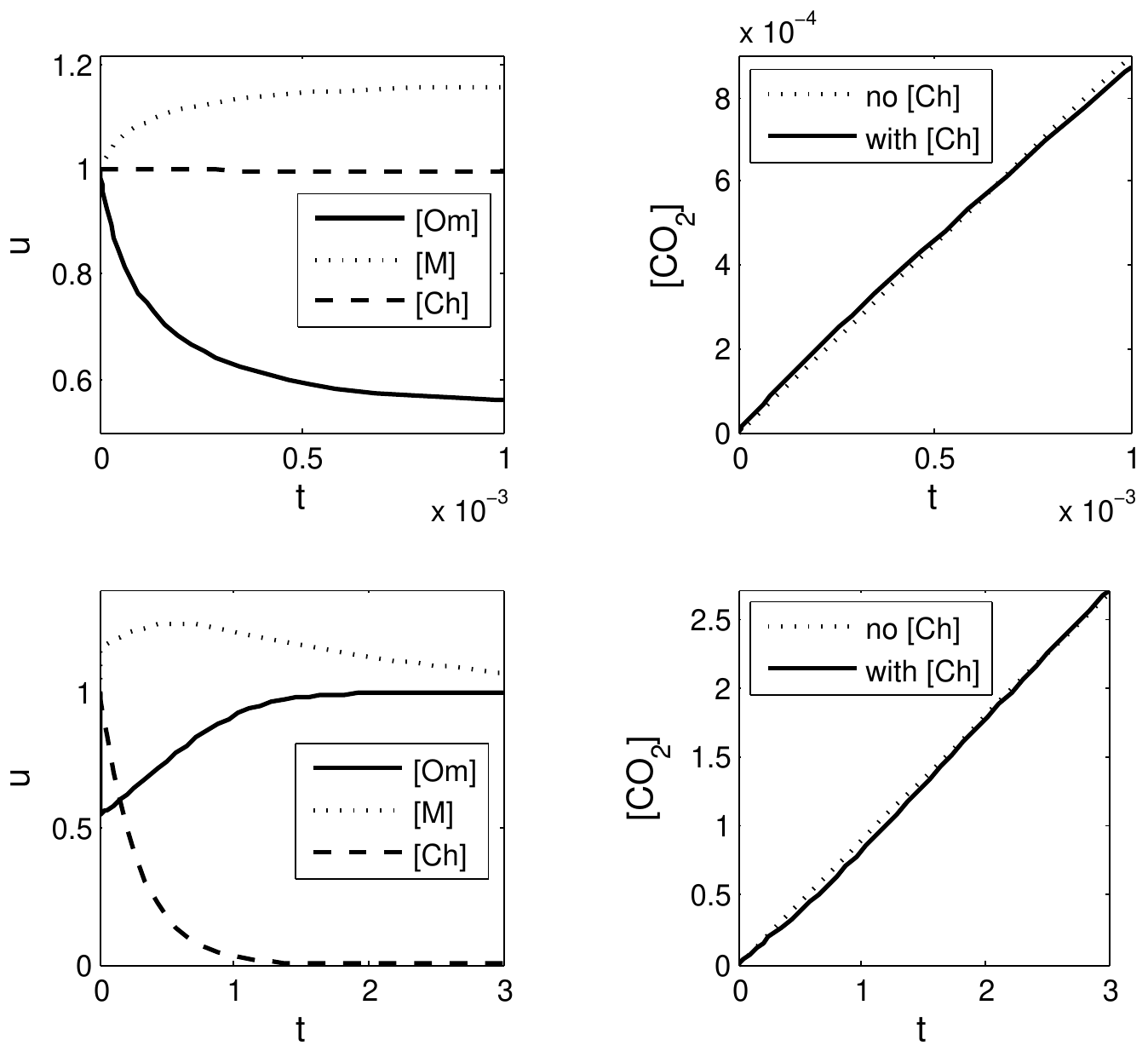}
                 \caption{\label{pmsetthree}
                 Result of the simulation by using Parameter Set 3 (see Table \ref{tab:params:3}).}
                 \end{centering}
\end{figure}

The results of the simulation are shown in Figure \ref{pmsetthree}. We see similar behavior as in Figure~\ref{pmsetone}. Hence we expect that different soils can allow for more a prominent charcoal influence on the overall $CO_2$ emission.


\section{Conclusion}\label{Discussion}

Within this framework, we translated the problem of  charcoal sequestration in soils in terms of the large-time asymptotics of classical solutions to a set of nonlinear differential equations describing a lumped chemistry between charcoal and chemical composition of soils.

Proving basic results (positivity and $L^\infty$-bounds on concentrations, well-posedness, stability of steady states, etc.), we point out the nice mathematical structure of the system.

What concerns the system's sensitivity with respect to varying parameters, we observe strong effects especially on intermediate time scales.
Most importantly,
for a rather large range of parameter values, our simulations clearly indicate that the short-time behaviour of our system can be significantly different from the long-time behaviour.
Therefore, when testing experimentally the effect of adding charcoal to soil on $CO_2$ emission, it is dangerous to make judgements based solely on short-time data.

Our model incorporates a large set of parameters and reference constants. To identify reasonable numerical ranges for them, we would need more experimental data (i.e. the $CO_2$ {\em vs.} time curve of other variations of $U_1,U_2$). A proper parameter identification work would naturally lead to a better control of the size of the characteristic time scales and potentially allow for improved predictions on $CO_2$ sequestration.

Our simulation output shows that  there is enough freedom to change the relative difference in $CO_2$ emission (with respect to not putting charcoal into the soil) both qualitatively and quantitatively. It seems however that more model components are needed to get better predictions. For instance,  the system of equations can be enlarged to include nutrients, phenols, temperature, etc.  and/or rain effects [maybe imposing some minimal spatial dependence in the model parameters, eventually also allowing for transport mechanisms].


\section*{Acknowledgments}
We thank the two reviewers for helping us improve our manuscript. Thanks extend also to C. Lazcano (Vigo) for posing us the Biochar problem, T. van Noorden (Gouda) and F. Wirth (W\"urzburg) for useful discussions, and  to C. van Altena (Wageningen)  for helping us with a better insight into soil data.  A.M. acknowledges support from RING (Research
Initial Network Grant) from British Council (France office) and a PPS RV22 award (Partnership Programme in Science) from  British Council (UK) and Platform B\`eta Techniek (NL).
P. v. M. is financially supported by the Complexity program of NWO (Netherlands Foundation for Scientific Research).


\appendix

\section{Physical meaning of  \eqref{dgem}: Equilibria and stability of the reaction block  \eqref{CR8} and \eqref{CR10}} \label{SAsub1}

We consider here the subsystem of \eqref{DE8}--\eqref{DE6} that corresponds to reactions
\eqref{CR8} and \eqref{CR10} (without the presence of charcoal, carbon dioxide or a source of organic matter).
The reason for studying this subsystem is that it gives us a physical reason for imposing \eqref{dgem}. Moreover, this subsystem turns out to dominate the short time behaviour of the whole system.

Substituting $s=0$, $k_1=0$ and $[Ch]=0$ into \eqref{DE8} and \eqref{DE9} gives
\begin{equation}\label{PR01}
\begin{aligned}
\frac{d}{dt} [Om] & = - \delta k_3(0) [M] [Om]^\delta + \eta k_4 [M], \\
\frac{d}{dt} [M] & =  \mu k_3(0) [M] [Om]^\delta - k_4 [M].
\end{aligned}
\end{equation}
In the rest of this subsection we write $k_3$ instead of $k_3(0)$ for brevity.

\begin{figure}[h!]
\begin{centering}
                 \begin{tabular}{lcr}
                  \includegraphics[width=150pt]{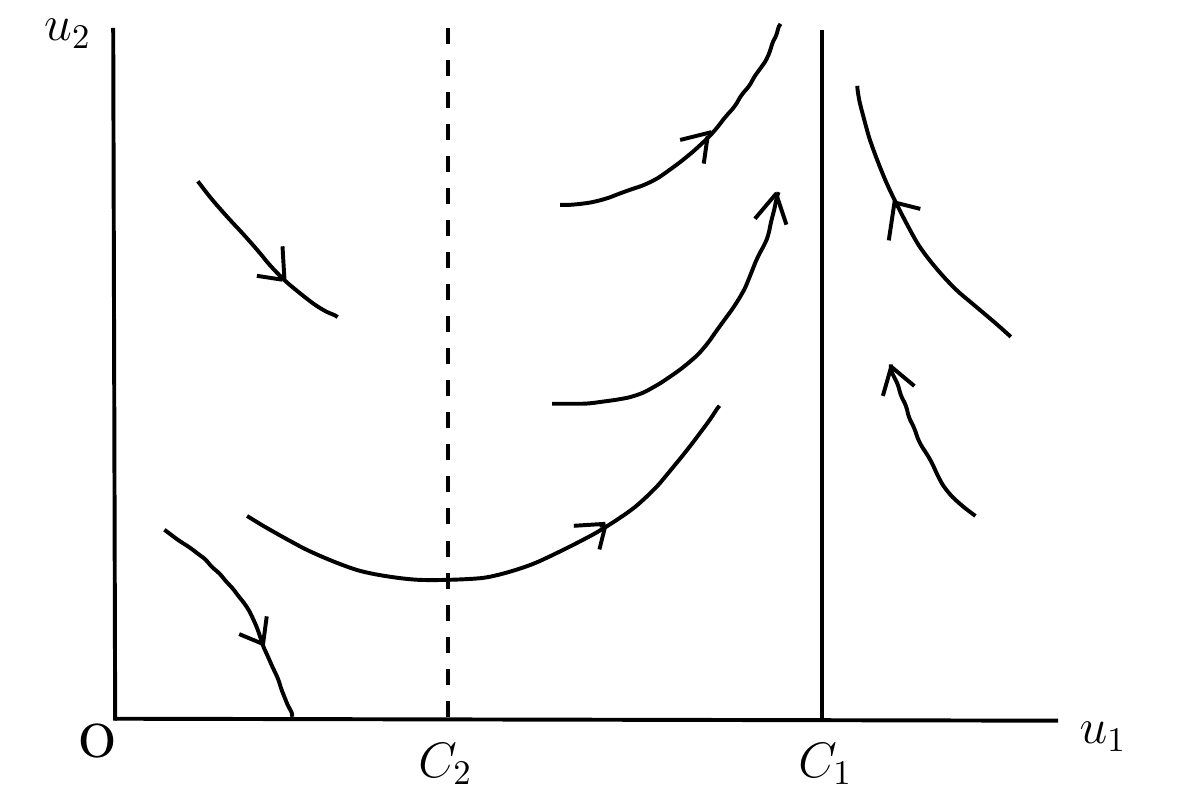} & \hspace{4mm} & \includegraphics[width=150pt]{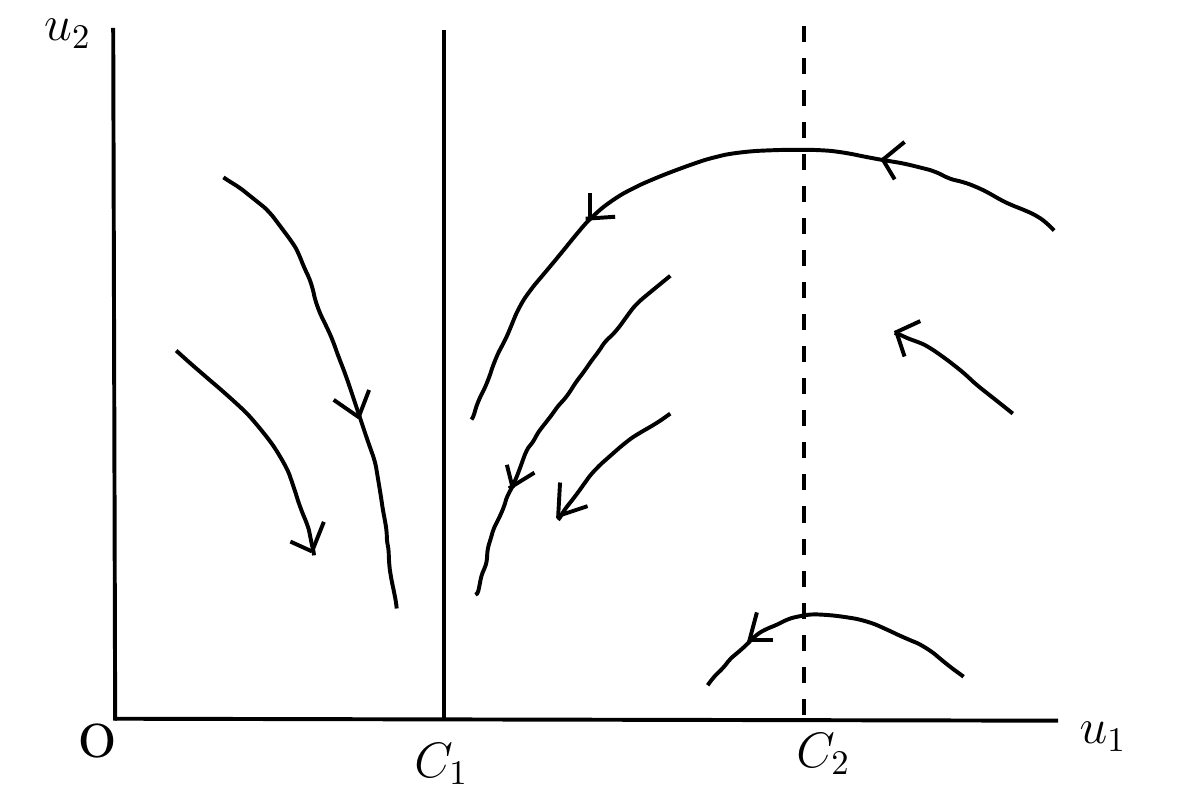} \hspace{8mm} \\
                 \end{tabular}
                 \caption{Sketches of the phase plane corresponding to \eqref{PR01}, depending on whether $C_1$ is bigger or smaller than $U_1$ (see \eqref{PR03} for their definitions). Recall that $u_1=[Om]$ and $u_2=[M]$.
                 }\label{fig1}
                 \end{centering}
\end{figure}

Figure \ref{fig1} shows a sketch of the phase plane corresponding to \eqref{PR01}.
Note that
\begin{equation}\label{PR02}
\begin{aligned}
    \frac{d}{dt} [Om] = 0 \hspace{1mm} &\Leftrightarrow \hspace{1mm} [M] = 0 \textrm{ or } [Om] =
    \left( \frac{ \eta k_4 }{ \delta k_3 } \right)^\frac{1}{\delta} =: C_1, \\
    \frac{d}{dt} [M] = 0 \hspace{1mm} &\Leftrightarrow \hspace{1mm} [M] = 0 \textrm{ or } [Om] =
    \left( \frac{ k_4 }{ \mu k_3 } \right)^\frac1\delta =: U_1.
\end{aligned}
\end{equation}
From \eqref{PR02} we see that $([Om],[M])=(c,0)$ is an equilibrium solution of \eqref{PR01} for all $c \in \mathbb{R}$. If $C_1=U_1$, then so is
$([Om],[M])=(C_1,c)$ for all $c \in \mathbb{R}$.

To determine the stability of the first equilibria, $([Om],[M])=(c,0)$,
we compute the Jacobian matrix
corresponding to the system \eqref{PR01}:
\begin{equation}\label{PR03}
   k_3 \begin{bmatrix}
     - \delta^2 [M] [Om]^{ \delta - 1 }       & \delta \left( \frac{\eta k_4}{\delta k_3} -
     [Om]^\delta \right) \\
     \delta \mu [M] [Om]^{ \delta - 1 } & \mu \left( [Om]^\delta - \frac{k_4}{\mu
     k_3} \right) \\
   \end{bmatrix}.
\end{equation}
From \eqref{PR03} it easily follows that the equilibria $([Om],[M])=(c,0)$ are stable if $c
< U_1$.

Now, we consider the boundedness of the trajectories.
We consider three cases: $C_1 < U_1$, $C_1 = U_1$ and
$C_1 > U_1$ (sketches of the corresponding phase planes are given in Fig. \ref{fig1}).
These cases correspond to:
\begin{description}
  \item[$(\delta > \eta \mu) :$] From the phase field analysis, we expect the solution of \eqref{PR01} to
      be bounded for all initial conditions.
  \item[$(\delta = \eta \mu) :$] From \eqref{PR02}, we see that we have more equilibrium
      points, which are given by $[Om] = C_1 = U_1$ and $[M] \in \mathbb{R}$ arbitrary. These
      equilibrium points are stably if and only if $[M] > 0$.
  \item[$(\delta < \eta \mu) :$] From the phase field analysis, we expect the solution to
      blow up for most initial conditions.
\end{description}
Therefore a sufficient condition for a solution of the reduced system \eqref{PR01} to be finite in time is
\begin{equation}\label{dgem2}
    \delta \geq \eta \mu.
\end{equation}
This is the same as our assumption \eqref{dgem} for the whole system.
Equality in \eqref{dgem2} would mean that the amount of organic matter that is converted into microbes by reaction \eqref{CR8} is equal to one over the amount of microbes that is converted into organic matter by reaction \eqref{CR10}. This means that $[Om] + \eta [M]$ is conserved. Indeed, one sees immediately from (\ref{PR01}) that
\begin{equation*}
    \frac{d}{dt} ([Om] + \eta [M]) = 0.
\end{equation*}
 This quantity $[Om] + \eta [M]$ was also useful for proving $L^\infty$ bounds for the whole system.
See equation \eqref{LB2}.

\section{Sensitivity of $CO_2$ emission for the Parameter Set 1}\label{sens}


Here we illustrate numerically how sensitive the behaviour of the $CO_2$ emission is with respect to changes in two values (one at a time) of Parameter Set 1. A similar discussion can be made based on Parameter Set 2 and Parameter Set 3. 

\begin{figure}[h!]
\begin{centering}
                  \includegraphics[width=350pt]{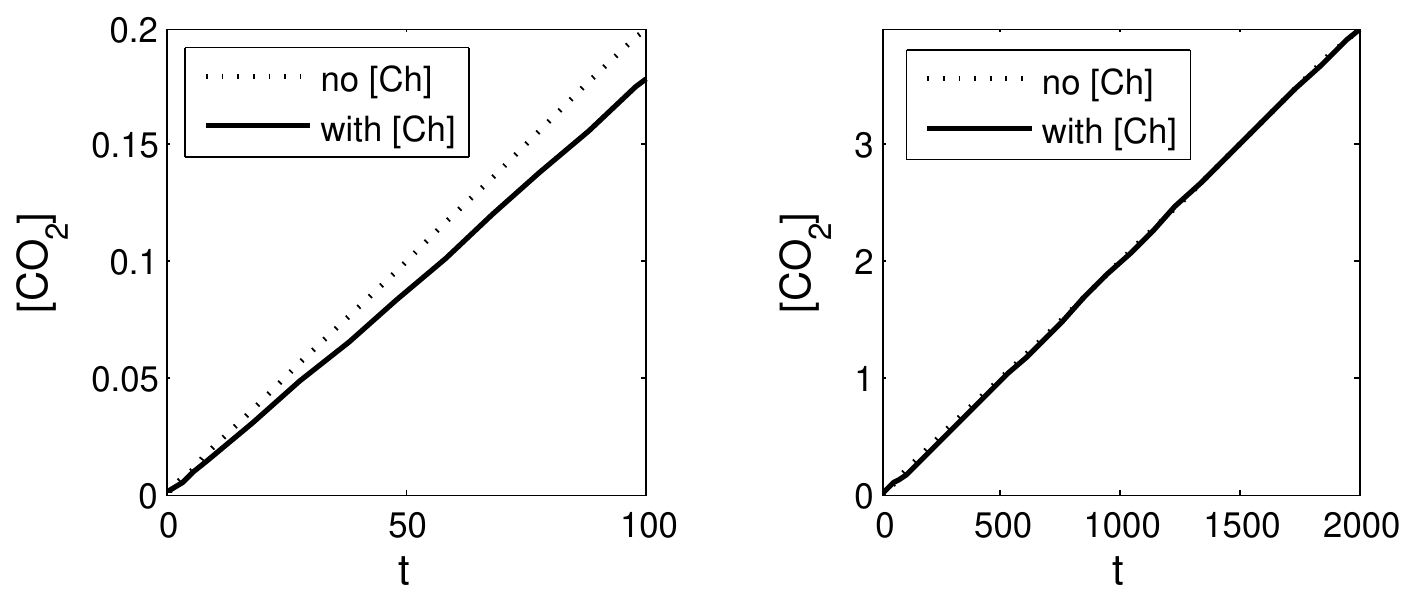}
                 \caption{ Parameter Set 1 with $U_3 = 10$ mol m$^{-3}$, i.e. ten times as much charcoal in the soil.}\label{fig6}
                 \end{centering}
\end{figure}

We start with testing the sensitivity of the $CO_2$ emission by increase the amount of charcoal that we put initially in the ground. We take $U_3 = 10$ mol m$^{-3}$ so that the amount of charcoal is ten times as much. The results are shown in Figure \ref{fig6}.

On a long time scale the behaviour is similar to before. This is remarkable, because it means that the total amount of emitted $CO_2$ hardly changes when ten times as much charcoal is put into the soil. On the short time scale we do see a difference: the rate of $CO_2$ emission is slightly increased, but it is still lower than the reference $CO_2$ emission.

\begin{figure}[h!]
\begin{centering}
                  \includegraphics[width=350pt]{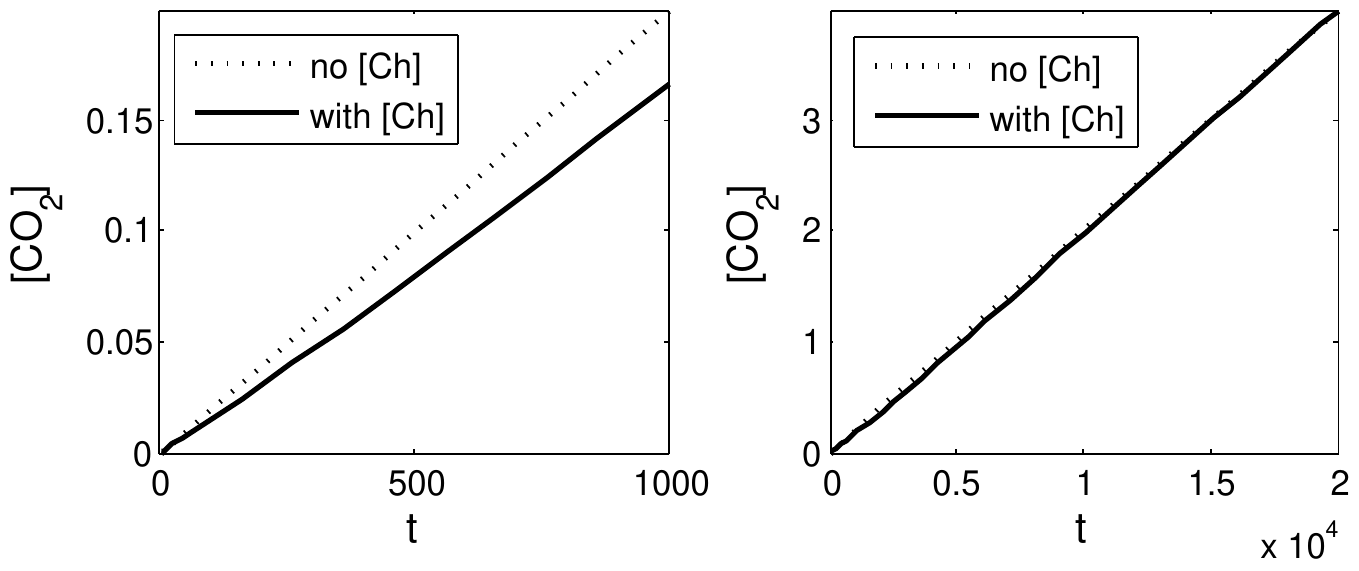}
                 \caption{\label{fig4}
                 Parameter Set 1 with $K_2 = 10^{-4}$ s$^{-1}$, i.e. the breakdown of the charcoal by the microbes is ten times as slow.}
                 \end{centering}
\end{figure}

Now we test the effect of $K_2$ on the $CO_2$ emission. We take $K_2 = 10^{-4}$ s$^{-1}$, which is ten times less as the value for $K_2$ in Parameterset 1. This corresponds to a slower breakdown of the charcoal by the microbes. Figure \ref{fig4} shows the results.

The qualitative behaviour of the $CO_2$ emission does not change, but now the corresponding time scales are larger. This can be explained by the charcoal being in the system for a longer time, which causes an increase in the rate at which equilibrium is reached.

\end{document}